\documentclass[12pt,leqno]{article}
\usepackage[cp1251]{inputenc}
\usepackage[T2A]{fontenc}
\usepackage{graphics}
\usepackage{color}
\everymath{\displaystyle}

\begin{document}

\begin{center}
\hspace{0.5cm} \textbf{\copyright 2009 A.M.Tsirlin}\\
\end{center}
\begin{center}
\noindent{\Large \bf Maximum Principle for variational problems with scalar argument\\}
 {\Large A.M.Tsirlin \\}
{Program  Systems  Institute,  Russian  Academy  of  Sciences, tsirlin@sarc.botik.ru} 
\end{center}

\subsection *{}
\small {~~~\textbf{Abstract}: In this paper the necessary conditions 
of optimality in the form of maximum principle 
are derived for a very general class of variational problems.
This class includes problems with any optimization criteria and constraints that 
can be constructed by combining some basic types (differential equation, 
integral equations, algebraic equation, differential equations with delays, etc).     
For each problem from this class the necessary optimality conditions are produced by 
constructing its Lagrange function $R$ and then by dividing its variables into three groups  
denoted as $u(t)$, $x(t)$ and $a$ correspondingly. $a$ are parameters which are constant over time. 
The conditions of optimality state that a non-zero vector function of Lagrange multipliers exists 
such that on the optimal solution function $R$ attains maximum on $u$, is stationary on $x$, and 
the integral of $R$ over the control period $S$ can't be improved locally. 
Similar conditions are also obtained for sliding regimes. Here solution is given by the limit of 
maximizing sequence on which the variables of the second group are switching with infinite frequency between 
some (basic) values.}

\section*{~~~1. Introduction}

The first formulation of the Pontryagin maximum principle was published in \cite{Pontr1} for a problem of optimal 
high-speed response. The maximum principle was then extended to various other classes of 
variational and optimal optimal control problems, see for example [2--5]. 
For each of these new classes the new proof of the maximum principle was derived. 
Most of these proofs used Rozonoer approach \cite{Rozon}, who published the first proper proof 
of the maximum principle based on using needle-like variations of control. 
The same proof was later used in Pontryagin's book  \cite{Pontr2}, where the  
needle-shaped variations were called the Makstein variations. 
Pontryagin did not mention neither Boltyanskii's \cite{Pontr1} nor Rozonoer's proofs, 
and did not acknowledge the crucial contribution of A.A. Fel'dbaum, who first stated the 
standard modern optimal control problem.

V.M. Tikhomirov repeatedly emphasized that it is desirable to consider various forms of 
optimal control problems within a unified framework by stating their optimality conditions 
in terms of the Lagrange function. This approach was implemented in \cite{Arut} 
for the problems described by the ordinary differential equations.

In practice, the optimal control problems may also have many other types of constraints 
- integral equations, finite constraints,  etc. The exact composition of these constraints 
may be very different in different parts of the control interval, etc.  

In this paper we will use the canonical form of a variational problem [8, 9] that includes 
a wide class of constraints and objective functions. We will prove the maximum principle 
for the problem in this form. The known proofs of the maximal principle rely heavily on 
the link between the problem's state (phase) and its control variables via the problem's 
differential constraints. The approach used in this paper is based on using the particular type of 
relaxation of the canonical extremal problem. Similar approach was applied by L.S. Yang, 
V.F. Krotov, V.I. Gurman, R.V. Gamkrelidze, and others [10--13] to obtaine bounds (estimates) 
on the solutions of the particular types of extremal problems, to investigate sliding regimes 
and to derive the sufficient conditions of optimality. 

When the particular problem from this class, with the particular type of constraints and the particular type of 
the optimality criterion, is considered, then this problem is  reduced it to the canonical form.  
This yields the problem's Lagrange function $R$ and defines which of its variables will belong to the first group and 
should maximize $R$ on the optimal solution. Thus, adding new constraint not only changes (adds a new term to) 
the problem's Lagrange function, but it also changes for which variables it is required to maximize it. 

\section*{2.Equivalent  transformations  and  relaxations  of  extremal  problems}

An extremal problem is defined as a problem of finding the maximum of
some criterion $I(y)$ on the feasible set $D$
\begin{equation}\label{98}
I(y) \to \max \bigg/ y \in D.
\end{equation}
$D$ can be a set in a vector space $R^n$ or in a space of
functions. The functional $I (y)$ is defined and bounded from above in $D$. 
The element $y^* \in D$, for which $I$ attains maximum
is called \textit{ the optimal solution}, and the value of
criterion $I(y^*)$ is called \textit{ the value of the problem}.
The set $D$ can belong to the vector space $R^n$ or to the space of functions. 
If the solution does not exist, then the value of the problem is the exact upper bound 
of the criterion $(\sup I (y))$ in $D$. In this case, the maximizing sequence is called the generalized solution.

\subsection*{2.1. Equivalent  transformations}

An extremal problem (\ref{98}) can be transformed into another extremal problem. 
If the solutions of both the original and transformed problems are the same then the
problems are called equivalent with respect to solutions. In their
values coincide then they are called equivalent with respect to
values. If both solutions and values are the same then the
problems are said to be equivalent. Let us give a few examples of such 
transformations.

(1) Applying monotone function to criterion transforms (\ref{98})
into the following problem
\begin{equation}\label{99}
F_0 \bigl( I(y) \bigr) \to \max \bigm/ y \in D,
\end{equation}
here $F_0$ is a monotonically increasing function (if $F_0$ is 
differentiable then its derivative on $I$ is positive almost
everywhere). The problem  (\ref{99}) is equivalent to (\ref{98})
with respect to the solution.

(2)  Adding vanishing term to the optimality criterion
\begin{equation}\label{910}
I(y)+ \varphi(y) \to \max \bigm/ y \in D,
\end{equation}
where the function $\varphi(y)$ is equal zero for any $y \in D$.
The problem (\ref{910}) is equivalent to the original one.

The concept of equivalent extremal problems can be generalized for
problems where feasible sets of the original and transformed
problems are different. Consider two extremal problems:

problem $A: \quad \; \; I_A(y) \to \max \bigm/ y \in D_A$;

and problem $A1: \quad I_{A1}(z) \to \max \bigm/ z \in D_{A1}$. 

Definition 1. \textit{ The problems $A$ and $A1$  are equivalent with
respect to solution if such one-to-one mapping between $D_A$ and
$D_{A1}$ can be found that from}
\begin{equation}\label{9.1}
I_A(y_1) \ge I_A(y_2), \quad (y_1, y_2) \in D_A \vspace*{-0.2cm}
\end{equation}
\textit{ it follows that }
\begin{equation}\label{9.2}
I_{A1}(z_1) \ge I_{A1}(z_2), \quad (z_1, z_2) \in D_{A1}.
\end{equation}
Here  $z_1$ is mapped onto $y_1$, and $z_2$ to $y_2$. We define
the class of equivalent problems $\overline A$ as a set of all
problem that are equivalent to $A$. The inequalities  (\ref{9.1}),
(\ref{9.2}) guarantee that the optimal solution of $A$ corresponds
to the optimal solution of $A1$.

\textit{ Example.} Suppose the problem $A$ and $A1$ have the
following forms
$$
f_0(y) \to \max \Bigm/ f_0(y) \geq 0, \quad a \le y \le b
$$
$$
\int\limits^b_a \sqrt{f_0(\tau)} \delta(t- \tau)d \tau \to \max
\Bigm/f_0(\tau) \geq 0.\vspace*{-0.2cm}
$$
here function $f_0$ is continuous and bounded on $[a, b]$. After
introducing mapping of solutions of the problem $A$ $y^0$ onto
solutions of the problem $A1$ $\delta(y^0- \tau)$ (where $\delta$
is the Dirac delta function), then $A1\in \overline A$. 

\subsection*{2.2. Relaxation of extremal problem}

The method we use in this paper to solve extremal problems 
includes transformation of the original problem into another problem 
with the same criterion and a larger feasible set [10, 11].  
The simplest transformation of this
type can be obtained by deleting one of the original problem's
constraints. It is clear that the new problem will have a wider
feasible set. The other transformations that can be used here include  
adding time-dependent points to the feasible set of the original problem 
that includes only time-independent points, 
allowing discontinuous solution in addition to smooth ones, etc. 
All these methods transform the original problem into its relaxation.

Since the original extremal problem $A$ corresponds to the entire class
$\overline A$ of the equivalent transformed problems, we shall call any
member of this class a relaxation of $A$ (unlike the conventional definition 
of relaxation \cite{Gurman}). 

Definition 1: \textit{The problem $B$ ($I_B(y) \to \max \bigm/ y
\in D_B$) is a relaxation of the problem $A$ ($I_A(y) \to \max
\bigm/ y \in D_A$), if it is possible to single out such subset
$\tilde D_B$ of $D_B$ that the problem $\tilde B$ ($I_B(y) \to
\max \bigm/ y $ $\in \tilde D_B$) is equivalent to the problem $A$
with respect to the solution}. Therefore the problem $\tilde B$
belongs to the class $\overline A$.

Note that feasible sets $D_A$ and $D_B$ can have different nature.
For example one can be a vector space and another a space of real
functions. In the general case $I_A$ and $I_B$ can be also different but
they obey the inequalities (\ref{9.1}), (\ref{9.2}) on $D_A$ and $\tilde
D_B$. In particular if $D_A$ coincides with $\tilde D_B$ then
$I_A$ and $I_B$ are either the same or one a is monotonic function
of another on $D_A$.

Sometimes criteria and constraints that determine $D$ depend on a
parameter $\lambda$ in such a way that for any $\lambda \in
V_{\lambda}$ the problem $B_{\lambda}$ is a relaxation of problem
$A$. We shall call such relaxation $B_{\lambda}$( $I_B(\lambda, y)
\to \max \bigm/ y \in D_{B \lambda})$  a \textit{ parametric
relaxation}.

Consider problem $A$: $I_A(y) \to \max \Big/ y \in D_A$ where a
finite set of constraints determines $D_A$. Suppose for each of
these constraints the norm $\Delta_j$ can be defined that measures its
deviation from some nominal value (that is, from its value for the original problem).

Definition 2: \textit{Problem $А$ is said to be well-posed if for
any $\epsilon >0$, $\delta$ exists such that from inequality
$\max_{j}(\Delta_j) \le \delta$ it follows that the absolute value of
the deviation of $I_A^*$ for problem with constraints
deviating from nominal values by $\delta_j$  from the problem's
value when all constraints have zero nominal values $\delta_j=0$ is
less than $\epsilon$.}

Definition 3: \textit{Relaxation  $B$: $I_B(y) \to \max \Big/ y
\in D_B$  of the problem $A$ is equivalent if }
\begin{equation}\label{4.1}
I^*_{\overline A}=\sup_{y \in D_{\overline A}} I_{\overline A}(y)=\sup_{y \in D_B}I_B(y)= I^*_{_B}.
\end{equation}

Note that the left hand side of this equality does not depend on
$I^*_{_A}$ and $D_{_A}$, it depends on $I^*_{_{\overline A}}$ and
$D_{_{\overline A}}$, that is, on the optimality criterion and
feasible set for any problem from the class of equivalent
relaxations of $A$.
 
\textit{Parametric relaxation is equivalent if equation }
(\ref{4.1}) \textit{holds for at least one $\lambda \in
V_{\lambda}$.}

The following statements follow from these definitions:

Lemma 1: \textit{The sufficient condition for relaxation to be
equivalent to the original problem is that for any solution of relaxation problem $y^0 \in
D_B$ it is possible to find a sequence $\{y_i\} \subset D_A$ of
the feasible solutions of the original problem such that}
\begin{equation}\label{4.2}
\lim_{i \to \infty} I_{\overline A}(y_i)=I_B(y^0).
\end{equation}
 	
For well-posed problems it is not necessary for $\{y_i\}$ to
belong to $D_{\overline A}$. It is only necessary that in the
limit $i\to\infty$ each of the constraints tends to nominal value
with arbitrary accuracy. Lemma 1 follows from the definition 3. If
optimal solution $y^*$ of relaxation problem exists then $y^0$ can
be replaced with $y^*$.

Lemma 2: \textit{If $y^*_A$ is the optimal solution of the problem
$A$,  $B$ is an equivalent relaxation of $A$ and $D_B \supset
D_A$, then the necessary conditions of optimality for relaxation
problem hold for $y^*_A$.}

Lemma 2 follows from the fact that $y^*_{A}$ can not be improved
on $D_A$, and $D_A$ is a subset of the feasible set of the
relaxation problem.

Relaxations are used

1) to reduce a conditional optimization problem to an
unconditional problem,

2) to find approximate solution in a class of maximizing sequences
if the problem does not have a solution in $D$,

3) to derive conditions of optimality and bounds on problem's
value,

4) to construct computational algorithms.

The most widely used type of relaxation uses criteria $I_A$ and
$I_B$ which have the same values on any element of $D_A$ and on
element that corresponds to it in  $\tilde D_B$. If $D_A$ and
$D_B$ are defined on the same space then  $D_A$ coincides with
$\tilde D_B$ and $\forall y\in D_A \quad I_A(y)=I_B(y)$. For this type
of relaxation the sufficient conditions of optimality (Krotov lemma \cite{KG}) holds.
It states that the sufficient condition for $y^*$ to be the solution of $A$ is that
$y^*$ is a solution of relaxed problem $B$ and $y^*$ belongs to $D_A$.

If the initial problem has no solution, then its is sufficient for 
a sequence of feasible solutions of the initial problem to 
be problem's generalized solution (maximizing sequence)
if this sequence approximates the solution of the relaxed problem with arbitrary accuracy.

\section*{3.Canonical form of the variational problem and optimality conditions for sliding regimes}

\subsection*{3.1. Variational  problem  in  the  canonical  form}

We shall call canonical the following variational problem:

\begin{equation}\label{4.3}
I=\int\limits^T_0 \Bigl[ f_{01}(t, x(t), u(t), a)+\sum\limits_l f_{02}(t, x(t), a) \delta (t-t_l)\Bigr] dt \to \max
\end{equation}                      
subject the  following conditions
\begin{equation}\label{4.4}
J_j(\tau)=\int\limits^T_0 \Bigl[ f_{j1}(t, x(t), u(t), a, \tau)+f_{j2}(t, x(t), a, \tau)\delta (t-\tau)\Bigr] dt=0,
\end{equation}
$$
\forall \tau \in [0, T], \quad j=1, \dots, m, \quad u
\in V_u,\quad a \in V_a,
$$
where $a$ is  vector of parameters that are constant on  $[0,
T]$; $u(t)$ and $x(t)$ are the piecewise continuous and the piecewise linear vector functions; 
values of $u(t)$ belong to the closed bounded set $V$ in the space $R^n$; 
functions $f_{j1}$ and $f_{j2}$, $j=0,\dots,m$ are defined on
direct product of the feasible sets of its two arguments, they are
both continuously differentiable on  $x$, $a$ and $t$ and $f_{j1}$ are continuous on $u$.
\\
$u(t)$ denote variables that are among the lists of parameters of 
the functions $f_{j1}$ for $j = 0, 1,..., m$ only (that is, $f_{j2}$ are independent on $u(t)$). 
We shall them the variables of the first group.
		
Lemma 3: \textit{Suppose the problem } (\ref{4.3}), (\ref{4.4})
\textit{is well-posed with respect to problem's value (according
to with definition 2, where the norm of deviation from the nominal
value of every constraint} (\ref{4.4}) \textit{is defined as
$\Delta_j= \max_{\tau}|J_j(\tau)|$), then the average relaxation of
this problem }
\begin{equation}\label{4.5}
\overline I=\int\limits^T_0\Bigl[\overline{f_{01}(t, x, u, a)}^u+\sum\limits_l
f_{02}(t, x, a) \delta (t-t_l)\Bigr] dt \to \max
\end{equation}
\textit{subject to}
\begin{equation}\label{4.6} \overline
J_j(\tau)=\int\limits^T_0 \Bigl[\overline{f_{j1}(t, x,u , a,
\tau)}^u+f_{j2}(t, x, a, \tau)\delta (t-\tau)\Bigr] dt=0
\end{equation}
$$
\forall \tau \in [0, T], \quad j=1, \dots, m, \quad u \in
V_\lambda, \quad a \in V_a
$$
\textit{ is equivalent to the problem (\ref{4.3}), (\ref{4.4})}.
\\
\\
\textit{Here}
\begin{equation}\label{4.7}
\overline{f_{j1}}^u=\int\limits_{V_u}
f_{j1}(t, x, u, a, \tau)P(u, t)du .
\end{equation}
The probability density measure $P(u,t)$ obeys the condition
\begin{equation}\label{4.8}
P(u, t)\ge 0, \quad \int\limits_{V_u} P(u, t)du=1 \quad \forall t \in
[0, T].
\end{equation}
This relaxation is obtained by replacing the control variable $u(t)$ with its mean value, that is, 
by applying randomization operation to $u(t)$. It was used by Yang [12] to study 
problems with differential constrains. It was also used by Krotov and Gamkrelidze to 
study sliding regimes in optimal control problems.

The  proof  of  the  Lemma  3  is  given  in  the  Appendix.

\subsection*{3.2. Optimality  condition  for  sliding  regimes}

The unknown variables in the problem (\ref{4.5})--(\ref{4.6}) is the measure $P^*(u, t)$, 
the vector function $x^* (t)$, and the vector $a^*$ . The solution of this problem 
obeys the following condition.
\\
\textbf{Theorem 1 (optimality  conditions  for  sliding  regimes for  the  variational problem  in  the  canonical  form)}.   
{\it Suppose}
$P^*(u, t), x^*(t), a^*$ {\it is the solution of the problem (\ref{4.5})--(\ref{4.6}).  Then  }
\\
1. {\it The optimal  distribution  of  the  randomized  variables  has  the following form}
\begin{equation}\label{4.9}
P^*(u,t)=\sum\limits^m_{\nu=0} \gamma_\nu(t)\delta(u-u^\nu(t)),
\end{equation}
\textit{where  for  all} $\forall t \in [0, T]$ 
\textit{the  piecewise  continuous  functions} $\gamma_\nu(t)$ \textit{obey  the  condition} $\gamma_\nu(t)\ge 0$, $\sum\limits^m_{\nu=0}
\gamma_\nu(t)=1$. 	
\\
2. \textit{We can find a scalar} $\lambda_0\ge 0$ \textit{and a continuous vector function} $\lambda(\tau)=(\lambda_1(\tau),...,\lambda_m(\tau))$, \textit{which are not equal zero simultaneously with} $\lambda_0$
\textit{on the interval} $[0, T]$ \textit{and are equal zero outside it, such that for the functional}
\begin{equation}\label{4.10}
S=\lambda_0 \overline{I} +\sum\limits^m_{j=1} \int\limits^T_0
\lambda_j(\tau)\overline{J_j}(\tau)d\tau=\int\limits^T_0 R dt
\end{equation}
\textit{and  its  integrand}
\begin{equation}\label{4.11}
R=\lambda_0 R_0+\sum\limits^m_{j=1} R_j,
\end{equation}
$$
R_0=\sum\limits^m_{\nu=0}\gamma_\nu(t)
f_{01}(t, x(t), u^\nu(t), a)+\sum\limits_l
f_{02}(t, x(t), a)\delta(t-t_l),
$$
\begin{equation}\label{4.12}
R_j=\int\limits^T_0 \lambda_j(\tau) \Bigg[
\sum\limits^m_{\nu=0}\gamma_\nu(t) f_{j1}(t, x(t), u^\nu, a,
\tau)+f_{j2}(t, x(t), a, \tau)\delta(\tau-t)\Bigg] d\tau
\end{equation}
\textit{the following conditions hold}
\begin{equation}\label{4.13}
\frac{\delta S}{\delta a}\delta a \le 0,
\end{equation}
\begin{equation}\label{4.14}
\frac{\delta R}{\delta x}=0,
\end{equation}
\begin{equation}\label{4.15}
u^\nu(t)=\arg\max_{u \in V_u}R(x, \lambda, a^*, u),\quad \nu=0,...,m,
\end{equation}
where $\delta a$ is  an admissible  variation  of  the  parameter $a$.
\\
$u^\nu(t)$ are called the {\it basic values} of the vector function $u$. 
Some of the conditions of the problem (\ref{4.5})--(\ref{4.6}) may not contain variables of the first group; in this case, 
the maximal number of basic values is less than $m+1$. The proof of the theorem is given in the Appendix.
\\
Note that the relaxation (\ref{4.5})--(\ref{4.6}) is equivalent to the original problem (\ref{4.3}),
(\ref{4.4}). Therefore from Lemma 2 it follows that if the optimal solution $(u^*(t), \linebreak x^*(t), a^*)$ of the initial problem 
exists in the class of piecewise continuous functions $u(t)$, then it satisfies the optimality conditions (\ref{4.13})--(\ref{4.15}). 
In this case $\gamma_0(t)=1$ and the remaining multipliers $\gamma_j(t)$ in (\ref{4.9}) are equal to zero.
 
The  optimality  conditions  of  the  problem  (\ref{4.3}),
(\ref{4.4})  takes  the  form  (\ref{4.13})--(\ref{4.14}) with 	
\begin{equation}\label{4.15i}
u^*(t)=\arg\max_{u \in V_u}R(x, \lambda, a^*, u),
\end{equation}
where	
\begin{equation}\label{4.11i}
R_0=f_{01}(t, x(t), u(t), a)+\sum\limits_l
f_{02}(t, x(t), a)\delta(t-t_l),
\end{equation}
\begin{equation}\label{4.12i}
R_j=\int\limits^T_0 \lambda_j(\tau) \Bigg[
f_{j1}(t, x(t), u(t), a, \tau)+f_{j2}(t, x(t), a, \tau)\delta(\tau-t)\Bigg] d\tau.
\end{equation} 	

Thus, for the variational problem in the canonical form (\ref{4.3}),
(\ref{4.4}) the maximum principle conditions (\ref{4.13}), (\ref{4.14}), (\ref{4.15i}) hold,
where the function $R$ includes contributions from the optimality citerion $R_0$ and from each of the constraints $R_j, j=1,2,...,m$.

\section*{4. Maximum principle for variational problems}

\subsection*{4.1. Necessary  conditions of optimality  for  sliding  regimes}

The theorem 1 ((\ref{4.13}), (\ref{4.14}), (\ref{4.15i})) allows us to
obtain the conditions of optimality in the form of maximum principle
for  problem with various types of criterion and constraints. 
This is done in two steps. Firstly, the problem under consideration is 
reduced to the canonical form as we it was done in the previous section, 
yielding $R_0$ term for criterion and $R_j$ terms for of its constraints. 
Then the variables of the first group with respect to every one of these terms
are singled out. The variables which belong to the first group with respect to all 
of these terms are denoted as $u(t)$. The maximum principle with respect to these variables (\ref{4.15i}) holds.


This process can be simplified if we derive contributions to $R_0$ and $R_j$ from the 
standard types of optimization constraints and optimality criteria together with  
the rules showing how to classify problem's variables for each such contribution. 
These contributions  and rules are shown in Tables 1 and 2.

We illustrate the derivation of such contributions to $R$ for problem with constraints in the form of the 
ordinary differential equation 
$$
x(\tau)=x_0+\int\limits_0^\tau f(x(t),u(t),t)dt.
$$
It can be rewritten using the $\delta$ function and the Heaviside step function $h(t)$ in the 
\newpage
{\footnotesize
Table  1. 
\vspace*{-0.2cm}
\begin{center}
\textbf{ Optimality  criteria  and  the  contibution to Lagrange function $R$ \\ from them  }\vspace*{-0.2cm}
\end{center}
\centerline{
\begin{tabular}{|l|l|l|l|}
\hline
№ & Optimality  criterion&  	Term $R_0$ & Type  of  \\
 &  $I\rightarrow \max$ &
 & term  \\
\hline
1  & $\int\limits_0^T f_0(y(t),a,t)dt$  &  $\lambda_0f_0(y(t),a,t)$ &
$R_{0I}$ \\
\hline
2  & $F_0(y(t_0),a,t_0)$ & $\lambda_0F_0(y(t),a,t)\delta(t-t_0)$ &
$R_{0II}$ \\
\hline
3 &I=$\min\limits_{t\in[0,T]}f_0(y(t),a,t)$ &
$\frac{\lambda_0I}{T}+\lambda(t)[I-f_0(y(t),a,t)]$ & $R_{0II}$ \\
   &   &  $\lambda(t)\le 0$, &  \\
   &   &  $\lambda(t)[I^*-f_0(y(t),a,t)]=0.$  &  \\
\hline
\end{tabular}
}
\vspace*{0.3cm}
Note:  If \quad  $I=\sum\limits_k a_k I_k$, then  $R_0=\sum\limits_k a_k R_{0k}$.}

{\footnotesize
Table  2.\vspace*{-0.2cm}
\begin{center}
\textbf{Basic  types  of  constraints and  their contributions to  $R$}\vspace*{-0.2cm}
\end{center}
\centerline{
\begin{tabular}{|l|l|l|l|}
\hline
№ &  Kind     &   & Type  \\
 & of  constraint & Term  $R_{j}$ &   of  term. \\
\hline
1  & $\int\limits_0^T f(y(t),a,t)dt=0$  &  $\lambda f(y(t),a,t)$ $\; \mbox{at} \;$  $t\in (0,T)$ & $R_{j}$ \\
      &         &   0 $\quad \mbox{at} \quad$ $t\notin (0,T)$   &    \\
\hline
2  & $f(y(t),a,t)=0,$ & $\lambda(t)f(y(t),a,t)$  $\; \mbox{at} \; $  $t\in(0,T)$ & $R_{j}$ \\
   & $\forall t\in(0,T)$  & 0 $\quad \mbox{at} \quad$ $t\notin (0,T)$ & \\
\hline
3 & $f(y(t_0),a,t_0)=0$ & $\lambda f(y(t),a,t)\delta (t-t_0)$ & $R_{j II}$ \\
\hline
4 & $\dot x=f(x(t),u(t),a,t)$ & $\psi(t)f(x(t),u(t),a,t)$, & $R_{j I}$ \\
   &    &  $\psi(t)=0 \; \mbox{at} \; t\notin [0,T]$  &  \\
  &  $\mbox{at} \; t\in[0,T]$  &  $\dot\psi(t)x(t)+(x(0)/T)\psi(0)$  & $R_{j II}$\\
\hline  
5 &$x(t)\!=\!\!\int\limits^t_0 \!f(x(\tau), u(\tau), \tau) d \tau + x_0$ &$f (x, u, t)\int\limits^t_T \lambda(\tau)d\tau $&$R_{j I}$\\
     &    &  $\lambda(t)=0 \; \mbox{at} \; t\notin [0,T]$  &  \\
     &    &$\lambda(t)(x(t)-x(0))$& $R_{j II}$\\
\hline
6 & $x(t)=\int\limits^T_0f(x(\tau),u(\tau),a,\tau,t)d\tau$  & $\int\limits^T_0\lambda(\tau)f(x(t),u(t),a,t,\tau)d\tau $  & $R_{сj I}$  \\
   &    &~~~~~~~~~ $-\lambda(t)x(t)$  & $R_{j II}$ \\
\hline
\end{tabular}}}
\newpage	
following form
$$
J(\tau)=\int\limits^{\overline t}_0\left[x(t)\delta(\tau-t)-
f(x(t),u(t),t)h(\tau-t)-\frac{x_0}{\overline t}\right]dt=0, \quad
\tau\in[0,\overline t].
$$
From (\ref{4.12}) it follows that  $R_j$ for this constraint has the form
$$
\begin{array}{l}
R_j=\int\limits^{\overline t}_0\lambda(\tau)\left[x(t)\delta(\tau-t)-
f(x,u,t)(h(\tau-t)-\frac{x_0}{\overline t})\right]d\tau= \\
=x(t)\lambda(t)-f(x,u,t)\int\limits^{\overline t}_t\lambda(\tau)d\tau-
\frac{x_0}{\overline t}\int\limits^{\overline t}_0\lambda(\tau)d\tau.
\end{array}
$$
Let us define the function $\psi(t)$ such that
$$
\dot\psi(t)=\lambda(t),\quad  \psi(\overline t)=0,\quad
\int\limits^{\overline t}_t\lambda(\tau)d\tau=-\psi(t), \quad
\int\limits^{\overline t}_0\lambda(\tau)d\tau=-\psi(0).
$$					
We get	
\begin{equation}\label{5.v64}
R_j=\dot\psi x+\psi f(x,u,t)+\psi(0)\frac{x_0}{\overline t}.
\end{equation}											
If the interval $[0,\overline t]$  and  $x_0$ are fixed, the last term in (\ref{5.v64})  has no effect on the optimality condition.

Other summands in the Tables 1 and 2 are obtained in the similar way. If necessity these tables can be 
extended for other types of constraints. 

Suppose we consider the problem of miximization of  one of the optimization criterion listed in Tabvle 1. 
subject to any combination of constraints from Table 2. 
We classify unknown variables of the problem as the variables of the first group 
is they are among the parameters of the functions $R_{0I}$ and $R_{0j}$ for every 
optimality criterion and every constraint. Thus the optimality conditions for the particular problem of maximization of the criterion $I$ 
shown in Table 1 subject to any combination of constraints from Table 2 can be obtained by  
summing up contributions to $R$ from the corresponding terms in these two tables, and then by 
dividing problem's uknown variables into two groups using the following rule: \textit{the problem's variable is classified as the one from the 
first group if $R_{0I}$ and $R_{j I}$ for the optimality criterion $I$ and for each of the constraints depend on this variable. 
All the other unknowns are classified as belonging to the second group.}

We will denote the variables of the first group by $u(t)$ and the variables of the second group by $x(t)$. 
The problem may have no variables of the first group if, for example, all the variables are linked to 
each other via the finite equation (row 3, Table 2).

Using  the  Tables  1  and  2  for  each  particular problem,  we  construct  the  function  $R$  as
\begin{equation}\label{5.v38}
R=\lambda_0R_0+\sum\limits^n_{\nu=1}R_{\nu}.
\end{equation}
We denote all terms in (\ref{5.v38}) that depend on $u(t)$ as $H$. 
We denote the rest of this expression as $N$. Thus,

$$R=N(x,\lambda,t)+H(x,u,\lambda,t).$$ 

Then the Lagrange  function  for  the  extended  problem becomes:

$$
\tilde R=N(x,\lambda,t)+\sum\limits^m_{k=0}\gamma_k(t)H(x,u_k,\lambda,t),
$$
where $m$ is the total number of the conditions of the problem, which contain $u(t)$, $m\le n$, while the conditions \vspace*{0.3cm}
$\gamma_k(t)$ satisfy the relations
\begin{equation}\label{5.v39}
\gamma_k(t)\ge 0, \quad \sum\limits^m_{k=0}\gamma_k(t)=1.
\end{equation}
	
\textit{The  optimal  solution  of  this  problem  in  the  class  of  sliding  regimes 
is defined as  such  functions $\gamma^*(t)$								
with  components $\gamma^*_k(t)$ and  $u^*(t)$  with  components  $ u^*_k(t)$,  and  also $ x^*(t)$, 
such that $u^*_k\in V, \: \gamma^*_k(t)$  satisfies	(\ref{5.v39}), while  the  vector  function $x^*(t)$ 	for  any $t\in[0,\overline t]$ 	can  be  approximated	as  exactly  as
desired by the sequence  $\{x_r(t)\}$ of admissible solutions (by constraint equations) of the 
problem such that in this solution the functional $I$ tends to its upper bound}. 
The measure of closeness of the functions $x^*(t)$  and $x_r(t)$   is the value that is maximum on $t$ 
of the absolute value of their difference.

The necessary  optimality  conditions  in  the  class  of  sliding  modes  are given by  the  following

\textbf{Statement} (corollary from the theorem 1). \textit{If $\gamma^*_k(t), \: u^*_k(t) \:
(k=0,\dots,m), \: x^*(t)$  is the solution of the problem of maximization of the functional $I$, listed 
in the Table 1, over the set of admissible solutions, which is defined by the conditions from 
the Table 2 in the class of sliding regimes, then there exists a vector function 
$\lambda=(\lambda_0,\lambda_1(t),\ldots,\lambda_m(t)); \: \lambda_0=(0;1)$, 
that is not equal zero for   $t\in[0,\overline t]$  and is equal zero outside of $[0,\overline t]$. 
For this vector function for almost all $t\in[0, \overline t]$,and $\gamma_k(t)\ge 0$ 
the function}
\begin{equation}\label{5.v40}
H(x^*,u^*_k,\lambda,t)=\max\limits_{u\in V}H(x^*,u,\lambda,t),
\end{equation}
\textit{attains its global maximum  on   $u$  at  $u^*_k$,  and the  extended  function 
$\tilde R$  is  stationary  on $x$}:
\begin{equation}\label{5.v41}
\frac{\partial N(x,\lambda,t)}{\partial x}=-
\sum\limits^m_{k=0}\gamma_k(t)\frac{\partial H(x,u_k,\lambda,t)}{\partial x}.
\end{equation}
\textit{If the problem also depends on the vector of parameters $a$, then the conditions of optimality 
include equations for finding the optimal value $a^*$ from the condition 
that the functional $S$ can not be locally improved with respect to $a$}:
\begin{equation}\label{5.v42}
\frac{\partial S}{\partial a}\delta a=\left[\frac{\partial}{\partial a}\int\limits^{\overline t}_0
\left[N(x,\lambda,a,t)+\sum\limits^m_{k=0}\gamma_k(t)H(x,u_k,a,\lambda,t)\right]
dt\right]\delta a \le 0,
\end{equation}
\textit{where $\delta a$ is  a feasible variation  of $a$}.				

All the functions that enter definitions of constraints of the problem must be continuous on $u$ and 
continuous differentiable on $x$ and $a$. 
For the existence of the maximum in (\ref{5.v40}), it is sufficient that the set $V$ be closed and bounded, 
and that the function $H$ be bounded by $u$.

\subsection*{4.2. Maximum  principle  for  problems  with  the  scalar  argument}

In $H$ has maximum on $u$ for almost all $t$ at a single point ($\gamma_0(t)\equiv 1$), then 
the problem has the solution $u^* (t)$ in the form of the piece-wise continuous function. 
From (\ref{5.v40}),(\ref{5.v42}) it follows that there exists a vector function 
$\lambda=(\lambda_0,\ldots,\lambda_m(t));$ $ \lambda_0=(0;1)$ that is not equal zero on 
$t\in[0,\overline t]$ and is equal zero outside of this interval. This vector function is such 
that at almost all $t\in[0, \overline t]$, the following conditions holds  
 \textbf{the maximum principle for problems with the scalar argument}:
\begin{equation}\label{5.v58}
H(x^*,u^*,t,\lambda,a^*)=\max\limits_{u\in V}H(x^*,u,t,\lambda,a^*),
\end{equation}
\begin{equation}\label{5.v59}
\frac{\partial N}{\partial x}=-\frac{\partial H}{\partial x},
\end{equation}
\begin{equation}\label{5.v60}
\frac{\partial S}{\partial a}\delta a=\left[\frac{\partial}{\partial a}\int\limits^{\overline t}_0
\left[N(x,\lambda,a,t)+H(x,u,a,\lambda,t)\right]
dt\right]\delta a \le 0.
\end{equation}
\subsection*{4.3. Examples}

\textbf{Pontryagin maximum principle.} 
As one of the examples of obtaining the necessary conditions 
of optimality in the form of (\ref{5.v58})---(\ref{5.v60}), we will consider the problem
\begin{equation}\label{5.v66}
\begin{array}{l}
I=\int\limits^{\overline t}_0f_0(x,u,t)dt+F_0(x(\overline t))\rightarrow \max,
\\
\dot x_\nu=f_\nu(x,u,t), \quad x_\nu(0)=x_{\nu 0}, \quad
\nu=1,\dots,m,\quad u\in V,
\end{array}
\end{equation} 
in which all constraints are in the form of differential equations. 

The Lagrange function here becomes 
\begin{equation}\label{5.v67}
R= \lambda_0R_0+ \sum\limits^m_{\nu=1}R_{\nu}= \lambda_0f_0+
\sum\limits^m_{\nu=1}(\psi_\nu f_\nu+\dot\psi_\nu
x_\nu)+\lambda_0F_0(x)\delta(t-\overline t).
\end{equation}
Note that the variables $u$ here enter only terms $R_{0II}$ and  $R_{j II}$. 
Hence, $u$ (controls) here belong to the variables of the first group. 
There are no parameters in the problem (\ref{5.v66}) 
and the optimality conditions (assuming that solution exists in the class 
of functions $x(t)$ that is differentiable for almost all $t$ and 
functions $u(t)$ that are piece wise continuous for almost all $t$) 
take the form
\begin{equation}\label{5.v68}
u^*(t)=\arg\max\limits_{u\in V}R(\lambda,u,x^*), \quad \frac{\partial
R}{\partial x_\nu}=0, \quad \nu=1,\dots,m,
\end{equation}
This  condition,  after taking into account (\ref{5.v67}),  can be reduced to the following conditions 		
\begin{equation}\label{5.v69}
\begin{array}{l}
u^*(t)=\arg\max\limits_{u\in V}H(\psi,u,x^*), \\
\dot\psi_\nu=-\frac{\partial H}{\partial x_\nu}=\lambda_0\frac{\partial
F_0}{\partial x_\nu}\delta(t-\overline t), \quad \nu=1,\dots,m.
\end{array}
\end{equation}
Here, the function $H$ (Hamilton function) is the sum of all terms in $R$ that depend on $u$:
$$
H=\lambda_0f_0+\sum\limits^m_{\nu=1}\psi_\nu f_\nu.
$$
	
If we take into account  that outside of the interval  $[0,\overline t]$  $\psi=0$, 
then from (\ref{5.v69}) it follows that $\psi_\nu(t)$ has a break at  $\overline t$ and 
\begin{equation}\label{5.v70}
\psi_\nu(\overline t)=\lambda_0\frac{\partial F_0}{\partial x_\nu}, \quad
\nu=1,\dots,m.
\end{equation}		
The conditions (\ref{5.v69}), (\ref{5.v70}) are to be solved together with the 
differential equations (\ref{5.v66}) 
and the boundary conditions for $x$. 

The proposed approach to obtaining the necessary conditions of optimality allows us to 
trace how these conditions  change when we add new constraints to the problem. 
For example, suppose we add to the problem (\ref{5.v66}) the condition  
\begin{equation}\label{5.v71}
F(x(\overline t))=0.
\end{equation}
This  adds the following term to  $R$  	
$$
\tilde R_{j}= \tilde \lambda F(x) \delta(t-\overline t).
$$
The  optimality  conditions  (\ref{5.v69}))do not change for $t<\overline t$.
For  $t=\overline t$ the  $\psi_\nu$ is now  equal  to						
\begin{equation}\label{5.v72}
\psi_\nu(\overline t)=\lambda_0\frac{\partial F_0}{\partial
x_\nu}+\tilde\lambda\frac{\partial F}{\partial x_\nu}, \quad
\nu=1,\dots,m.
\end{equation}	
The  additional  variable  $\tilde\lambda$ is  to be found from 
the  condition (\ref{5.v71}).

$F_0$ in  (\ref{5.v66}) and $F$ in (\ref{5.v71}) can also depend on controls. 
In this case the conditions of optimality do not change for $t<\overline t$. 
But the control $u(\overline t)$ now obeys weaker condition of local optimality 
$\frac{\partial}{\partial u}[\lambda_0F_0+\tilde\lambda F]\delta u\le 0,
$
because  it turned out that control here belongs to the variables of the second type for 
$\overline t$.
\\  
\textbf{Butkovskii optimality conditions for the problem with constraints in the form of 
integral equations}\cite{Butkovskiy}. For the problem
\begin{equation}\label{5.v73}
I= \int\limits^{\overline t}_0f_0(x,u,t)dt \rightarrow \max\limits_{u\in V_u}
\left/
\begin{array}{l}
\int\limits^{\overline t}_0f(x(\tau),u(\tau),t,\tau)dt-x(t)=0
\end{array}
\right.
\end{equation}
the  Lagrange   function  (using Tables  1  and  2)  has  the  form	
$$
R=\lambda_0f_0(x,u,t)+\int\limits^{\overline
t}_0\lambda(\tau)f(x(t),u(t),\tau,t)d\tau-\lambda(t)x(t).
$$
Here $u(t)$ is the variable of the first  group.  
The  optimality  conditions  will  take  the  form	
\begin{equation}\label{5.v74}
\begin{array}{l}
\displaystyle{\frac{\partial R}{\partial x}}=0\Rightarrow\lambda(t)=\frac{\partial}{\partial x}
\left[\lambda_0f_0+\int\limits^{\overline
t}_0\lambda(\tau)f(x,u,\tau,t)d\tau\right], \\
\\
u^*=\arg\max\limits_{u\in V}\left[\lambda_0f_0+ \displaystyle{\int\limits^{\overline t}_0} \lambda(\tau)f(x,u,\tau,t)d\tau\right].
\end{array}
\end{equation}

\textbf{Combination of differential and integral equations.} 
In many cases, it is more convenient to describe a linear object by a convolution equation rather than 
by a differential equation. We will show how the replacement of the differential equation 
(for definiteness, the equation for the $m$-th differential constraint) by the convolution equation 
of the form
$$
\int\limits^{\overline t}_0u(\tau)k(t-\tau)d\tau-x(t)=0.
$$
will effect the  necessary conditions of optimality (the maximum  principle).  
Here,  $k(t)$  is  an  impulse  transfer function.

Instead of the term 
$u_m(t)\int\limits^{\overline t}_0\lambda(\tau)k(\tau-t)d\tau$ the function $H$  
now includes the term $\psi_m f_m$. Similarly  the function $N$ now includes 
$(-\lambda_m x_m)$ instead of $(\dot\psi_m x_m)$. 
The division of the problem's variables between the first and the second groups do not change. 
The maximum principle here follow directly follow from (\ref{5.v58})--(\ref{5.v60}) and do not 
require any special derivation. 

\textbf{Problems with the conditions in the form of inequalities and with the maximin criterion.} 
Some of the conditions of a problem may have the form of inequalities. 
To obtain the optimality conditions using the proposed approach 
these inequalities can be rewritten in the form of equalities using additional artificial variables. 
For example, the inequality
\begin{equation}\label{5.v76}
f(y(t),t)\ge 0
\end{equation}
can  be  rewritten  as  the  equality
$$
f(y(t),t)-z(t)=0,
$$
using additional artificial nonnegative  variable  z(t).  
The  appropriate  term in  $R$  has  the  form
$$
R_{\nu}=\lambda(t)f(y,t)-\lambda(t)z(t).
$$

The variable $z(t)$ belongs to the second group and does not enter into other terms in $R$, 
except $R_{\nu}$. The conditions of local optimality of $R$ by $z$ (after taking into account 
that $Z$ is nonnegative and its feasible variation $\delta z\ge 0$) yields
$$
\frac{\partial R}{\partial z}\delta z\ge 0\Rightarrow\frac{\partial
R_{\nu}}{\partial z}\ge 0\Rightarrow\lambda(t)\ge 0.
$$
Here $\lambda(t)=0$  if $z(t)  >  0$,  that is,  if $f(y,t)>0$,  and 
$\lambda(t)>0$  if  $f (y,t)  =  0$. This is an exact analog of the 
condition of complementary slackness in  the mathematical  programming.	

For the maximin problem 
\begin{equation}\label{5.v77}
I=\min\limits_{t\in[0,\overline t]}f_0(y(t),t)\rightarrow \max.
\end{equation}
we can use the same method by adding an additional parameter $a$ to the problem which is  	
independent  of  $t$. Then the problem can  be  rewritten  as	
\begin{equation}\label{5.v78}
a\rightarrow \max
\end{equation}
subject to inequality that holds for any $t\in[0,\overline t]$,
\begin{equation}\label{5.v79}
f_0(y(t),t)-a\ge 0.
\end{equation} 

The criterion (\ref{5.v78}) and   condition (\ref{5.v79}) contribute  the following terms to the 
function  $R$  
$$
\tilde R=\lambda_0\frac{a}{\overline t}+\lambda(t)f_0-\lambda(t)a,
$$
where (similar to (\ref{5.v76}) $\lambda(t)\ge 0$, and
$$
\lambda(t)[f_0(t,y^*(t))-a^*]=0.
$$
Here, $\lambda_0=1$ if a non-degenerate solution exists. Otherwise $\lambda_0=0$.   
Because all terms in $R$ (except from $\tilde R$) do not depend on $a$ and because 
$a$ is unconstrainted, the condition of stationarity of the Lagrange functional $S$ on $a$ 
yields
\begin{equation}\label{5.v80}
\frac{\partial S}{\partial a}=\int\limits^{\overline t}_0\frac{\partial \tilde
R}{\partial a}dt=0\Rightarrow \int\limits^{\overline
t}_0\lambda(t)dt=\lambda_0.
\end{equation}

{\footnotesize
\begin{flushright}
{\bf APPENDIX}
\vspace*{-0.2cm}
\end{flushright}

\textbf{Validity of the Lemma 3} follows from the Lemma 1 and the fact that for 
any solution $P^0(u, t), \, x^0(t), \, a^0$ we can find the sequence of points 
$$
\{ z_i \} = \{ u_i(t), x^0(t), a^0 \},
$$
on  which
$$
I(z_i) \to \overline{I^*}, \quad J_j(\tau, z_i) \to \overline{J_j}(\tau)=0,
\quad i \to \infty, j= \overline{1,m}.
$$
Indeed, for any given  $t, x, a,
\tau$, the vector $\overline f= (\overline f_0, \overline f_1, ...,
\overline f_m)$ belong to the convex hull $\overline{Q}$ of the set $Q$ which is obtained 
by  mapping of $V_u$ onto $(m+1)$-dimensional space $f$.
The solution maximizes $f_0$ with respect to $u$. 
Therefore, it belongs to the upper bound $\overline{Q}$ and can be obtained as a 
linear combination of no more than $(m + 1)$ elements of $Q$ (Carathйodory's theorem).

For any solution $P^0(u, t)$ that has the form (\ref{4.8}), it is possible to 
construct the sequence  $\{ u_i(t) \}$ of solutions of the problem (\ref{4.3}), (\ref{4.4}) 
by dividing the interval $[0, \tau]$ into $i$ subintervals $\Delta_1, 
..., \Delta_i$ and assuming that $\gamma_{\nu}(t)$ and
$u^{\nu}(t)$ are constant on each of these intervals. 
Suppose their values are denoted as $\gamma_{\nu r}$ and $u^{\nu}_r \, (\nu = 0, \dots, m)$ 
correspondingly. We shall call the problem obtained in this way 
the discretization of the problem (\ref{4.5}), (\ref{4.6}).

We divide the interval $[0, \tau]$ in the problem (\ref{4.3}), (\ref{4.4}) in a similar way.
But here we divide each subinterval into $(m+1)$ smaller intervals. 
Thus, we divide $\Delta_r$ into $\Delta_{r0}, \Delta_{r1}, ..., \Delta_{rm}$, 
and we get $\Delta_{r \nu}/\Delta_{r}= \gamma_{\nu r}$. We assume that the variables $u(t)$ 
in the problem (\ref{4.3}), (\ref{4.4}) are piece-wise constant and are equat to $u^{\nu}_r$ 
on the interval $\Delta_{r \nu}$. 
For the solutions constructed in this way on each of the intervals $\Delta_{r}(r=1,...,
i)$, the values of the functionals $I$ and $J(\tau)$ in the problem (\ref{4.3}), (\ref{4.4}) 
are equal to the values of the corresponding functionals for the discretization problem 
(\ref{4.5}), (\ref{4.6}). If $i \to \infty$ then $\Delta_r$ tend to zero uniformly in $r$ 
and $I_D$ and $J_D(\tau)$  for the discretization of the averaged
problem become arbitrary close $\overline I (\tau)$ and $\overline J (\tau)$. 
Because the problem under consideration in well-posed (definition 2)  
the same is also true for $I$ and $J(\tau)$. The Lemma 3 is proved.

\textbf{For  the  proof  of  the  theorem}  we  use the following  statement:

\textit{Suppose $y^*(t)$  is  the  solution  of  the  following problem  }
\begin{equation}\label{p7}
I= \int\limits^T_0 f_0(y, t)dt \rightarrow \max,
\end{equation}
\textit{subject to  constraints}

\begin{equation}\label{p8}
J_j(\tau)= \int\limits^T_0 f_j(y, t, \tau)dt=0, \quad j= 1, \dots, m, \quad
\tau \in [0, T],
\end{equation}
\textit{where $f$ is continuous and continuously differentiable with respect to all arguments. 
Then a non-zero vector }
$$
\lambda=(\lambda_0, \lambda_1(\tau), ..., \lambda_m(\tau)), \quad \lambda_0
\ge 0,
$$
\textit{can be found such that  for $y  = y^*$  the following inequality  holds }		
\begin{equation}\label{p.8}
\left( \frac{\partial R}{\partial y} \right) \delta y \le 0,
\end{equation}
\textit{where}
$$
R=R_0+ \sum\limits^m_{j=1} R_j= \lambda_0 f_0+ \sum\limits^m_{j=1}
\int\limits^T_0 \lambda_j(\tau) f_j(y, t, \tau)d \tau;
$$
\textit{and where $\delta y$
 is  the feasible  variation  of  $y(t)$  with respect to the  condition $y \in V_y(t)$.}

Proof.~ For simplicity we assume that $m=1$. After expanding $f_0$ and $f_1$ near $y^*(t)$ and 
neglecting higher than linear terms we get:
$$
f_0(y,t)=f_0(y^*,t)+(\partial f_0/\partial y)\delta y,\qquad f_1(Y^*,t)+(\partial f_1 /\partial y)\delta y.
$$ 

Suppose the problem is non degenerate. That is, that $y^*(t)$ is not an extremal of $J_1(\tau)$ for 
any $\tau\in[0,T]$. Then such value $t=t_1(\tau)$ can be found that  
\begin{equation}
(\partial f_1(y^*,t_1,\tau)\partial y) \neq 0.
\label{p4}
\end{equation}

Consider the variation of the solution $\delta y$ which differs from zero only over two 
infinitizemal time intervals $\epsilon_t$ around $t_1(\tau)$ and the arbitrary $t_2 \in [0, T].$ We denote  
$$
\sigma_1=\int_{\epsilon_{t_{1}}}\delta y(t)dt,\qquad \sigma_2=\int_{\epsilon_{t_{2}}}\delta y(t)dt.
$$
Variations of $I$ and $J_1(\tau)$ take the form 
$$
\delta I=(\partial f_0/\partial y)_{t_{1}} \sigma_1 + (\partial f_0/\partial y)_{t_{2}}\ \sigma_2, \quad \delta J_1(\tau)=(\partial f-1/\partial y)_{t_{1}}\sigma_1+(\partial f_1/\partial y)_{t_{2}}\sigma_2.
$$
Since the latter expression is equal zero for arbitrary $\tau \in [0, T]$, we can rewrite it as 
\begin{equation}
\int\limits_0^T G(\tau)\delta J_1(\tau)d\tau=\sigma_1 \int\limits_0^T(\partial T_1/\partial y)_{t_{1}}G(\tau)d\tau+ \sigma_2\int\limits_0^T(\partial f_0/\partial y)_{t_{2}}G(\tau)d\tau=0
\label{p5}
\end{equation}
for arbitrary function $G(\tau)$, which obeys the conditions $G(\tau)\geq 0,~ \int\limits_0^T G(\tau)d\tau=1,$ and the variation of the functional $I$ over the set of variations $\delta y$, feasible with respect of conditions (\ref{p5}), must non-positive:
$$
\begin{array}{c}
\delta I=\sigma_2\Bigg[(\partial f_0/\partial y)_{t_{2}}-(\partial f_0/\partial y)_{t_{1}} \int\limits_0^T(\partial f_1/\partial y)_{t_{2}}G(\tau)d\tau \Bigg/ \\
\Bigg/ \int\limits_0^T(\partial f_1/\partial y)_{t_{1}}G(\tau)d\tau \Bigg] \leq0.\\
\end{array}
$$
From (\ref{p4}) it follows that the denominator of the fraction in the square bracket is non zero. Denoting  
$$
\lambda(\tau)=-\left(\frac{\partial f_0}{\partial y}\right)_{t_{1}}G(\tau) \Bigg/ \int\limits_0^T \left(\frac{\partial f_1}{\partial y}\right)_{t_{1}}G(\tau)d\tau 
$$
and taking into account that $t_2$ can have an arbitrary value, we get 
\begin{equation}
\frac{\partial}{\partial y}\left[f_0(y,t)+\int\limits_0^T\lambda(\tau)f_1(y,t,\tau)d\tau\right]\delta y\leq 0.   
\label{p6}
\end{equation}
If the condition (\ref{p4}) does not hold then a non-zero function $\lambda_1(\tau)$ can be found such that  
\begin{equation}
\frac{\partial}{\partial y}\left[\int\limits_0^T\lambda_1(\tau)f_1(y,t,\tau)d\tau\right] =0.
\label{p.7}
\end{equation}
Combining  (\ref{p6}) and (\ref{p.7}) yields the condition of optimality  (\ref{p.8}), where for a non-degenerate problem we can set $\lambda_0=1$. 	

The problem (\ref{4.5}), (\ref{4.6}) for the distribution $P(u, t)$ in the form (\ref{4.9}) has the form (\ref{p7}), (\ref{p8}), with
$$
\overline R= R-R_{m+1}= \lambda_0 R_0 + \sum\limits_j R_j - R_{m+1},
$$
where $R_0$ and $R_j$  have  the  form (\ref{4.12}),  and  the  term 
$R_{m+1}$  corresponds  to  the  condition
$$\sum\limits^m_{\nu=0} \gamma_{\nu}(t)-1=0 \quad \forall t \in [0, T],$$
which  can  be  rewritten  in  the  form (\ref{p8})  as
$$
J_{m+1}(\tau)= \int\limits^T_0 \left( \sum\limits^m_{\nu=0} \gamma_{\nu}(t)-1
\right) \delta(t- \tau)d \tau=0 \quad \forall \tau \in [0, T].
$$

Thus,
$$
R_{m+1}= \int\limits^T_0 \lambda_{m+1}(\tau) \left( \sum\limits^m_{\nu=0}
\gamma_{\nu}(t)-1 \right) \delta (t- \tau) d \tau= \lambda_{m+1}(t) \left(
\sum\limits^m_{\nu=0} \gamma_{\nu}(t)-1 \right).
$$

From the  conditions  (\ref{p.8})  for  $\gamma_{\nu}$
$$
\frac{\partial \overline R}{\partial \gamma_{\nu}} \delta \gamma_{\nu} \le 0,\quad
\quad \gamma_{\nu} \ge 0
$$
it follows that for  the basic  values $u^{\nu}(t)$ (where $\gamma_{\nu}(t) > 0$)
$$
R(x, \lambda, a^*, u^{\nu})= \lambda_{m+1}(t), \quad \nu= 0, \dots, m,
$$
and  for $u \ne u^{\nu}(t)$,  $\gamma_{\nu}(t)=0$  and $\delta \gamma_{\nu} > 0$,  and,  hence, 
$R(x, \lambda, a^*, u) \le \lambda_{m+1}(t)$. Therefore the 
maximum  condition  (\ref{4.15}) holds.

The condition (\ref{4.14}) follows from (\ref{p.8}) if we take into account 
that $x$ is unconstrained. The conditions (\ref{4.13}) follow 
from the fact that with respect to the vector of parameters $a$, the problem 
(\ref{4.5}), (\ref{4.6})is a nonlinear programming problem and $S$ is its Lagrange function.

\end{document}